\documentclass[11pt,a4,twoside,reqno]{amsart}

\usepackage{amsmath}
\usepackage{amssymb}
\usepackage{bm}
\usepackage{mathrsfs}
\usepackage{stmaryrd}

\usepackage{color}

\usepackage{url}

\usepackage{supertabular}

\newtheorem{thm}{Theorem}

\newtheorem{lemma}[thm]{Lemma}
\newtheorem{prop}[thm]{Proposition}

\newtheorem{problem}[thm]{Problem}

\numberwithin{equation}{section}

\DeclareMathAlphabet{\mathsfsl}{OT1}{cmss}{m}{sl}

\newcommand{\titl}{\textsl}
\newcommand{\term}{\emph}

\newcommand{\cnst}[1]{\mathrm{#1}}

	{\makebox{\phantom{}} \\ %
		\textsc{Input:}\begin{itemize}}%
	{\end{itemize}}

	{\textsc{Output:}\begin{itemize}}%
	{\end{itemize}}

	{\textsc{Procedure:}\begin{enumerate}}%
	{\end{enumerate}}

\renewcommand{\phi}{\varphi}

\newcommand{\eps}{\varepsilon}

\newcommand{\defby}{\overset{\mathrm{\scriptscriptstyle{def}}}{=}}

\newcommand{\econst}{\mathrm{e}}
\newcommand{\iunit}{\mathrm{i}}

\newcommand{\onevct}{\mathbf{e}}

\newcommand{\oneton}[1]{\left\llbracket {#1} \right\rrbracket}

\newcommand{\abs}[1]{\left\vert {#1} \right\vert}
\newcommand{\abssq}[1]{{\abs{#1}}^2}

\newcommand{\Prob}[1]{\operatorname{\mathbb{P}}\left\{ {#1} \right\}}
\newcommand{\Expect}{\operatorname{\mathbb{E}}}

\newcommand{\vct}[1]{\bm{#1}}
\newcommand{\mtx}[1]{\bm{#1}}

\newcommand{\adj}{*}

\newcommand{\rank}{\operatorname{rank}}

\newcommand{\diag}{\operatorname{diag}}
\newcommand{\trace}{\operatorname{trace}}

\newcommand{\norm}[1]{\left\Vert {#1} \right\Vert}

\newcommand{\smnorm}[2]{{\bigl\Vert {#2} \bigr\Vert}_{#1}}

\newcommand{\enorm}[1]{\norm{#1}_2}
\newcommand{\enormsq}[1]{\enorm{#1}^2}

\newcommand{\pnorm}[2]{\norm{#2}_{#1}}

\newcommand{\bigO}{{\rm O}}

\begin{document}

\title[Paving Uniformly Bounded Matrices]
{The random paving property for \\
uniformly bounded matrices}

\author{Joel A.\ Tropp}

\thanks{JAT is with Applied \& Computational Mathematics, MC 217-50, California Institute of Technology, 1200 E.~California Blvd., Pasadena, CA 91125-5000.  E-mail:
  \url{jtropp@acm.caltech.edu}.  This work was supported by NSF DMS 0503299.}
 
\date{14 December 2006.  Revised 29 January 2007, 4 September 2007, and 24 December 2007.  Accepted for publication, \titl{Studia Mathematica}.}

\begin{abstract}
This note presents a new proof of an important result due to Bourgain and Tzafriri that provides a partial solution to the Kadison--Singer problem.  The result shows that every unit-norm matrix whose entries are relatively small in comparison with its dimension can be paved by a partition of constant size.  That is, the coordinates can be partitioned into a constant number of blocks so that the restriction of the matrix to each block of coordinates has norm less than one half.  The original proof of Bourgain and Tzafriri involves a long, delicate calculation.  The new proof relies on the systematic use of symmetrization and (noncommutative) Khintchine inequalities to estimate the norms of some random matrices.
\end{abstract}



\maketitle

\section{Introduction}

This note presents a new proof of a result about the paving problem for matrices.  Suppose that $\mtx{A}$ is an $n \times n$ matrix.  We say that $\mtx{A}$ has an \term{$(m, \eps)$-paving} if there exists a partition of the set $\{1, 2, \dots, n\}$ into $m$ blocks $\{\sigma_1, \sigma_2, \dots, \sigma_m\}$ so that
$$
\norm{ \sum\nolimits_{j=1}^m \mtx{P}_{\sigma_j} \mtx{A} \mtx{P}_{\sigma_j} } \leq \eps \norm{ \mtx{A} }
$$
where $\mtx{P}_{\sigma_j}$ denotes the diagonal projector onto the coordinates listed in $\sigma_j$.  Since every projector in this note is diagonal, we omit the qualification from here onward.  As usual, $\norm{ \cdot }$ is the norm on linear operators mapping $\ell_2^n$ to itself.  

The fundamental question concerns the paving of matrices with a zero diagonal (i.e., \term{hollow matrices}). 

\begin{problem}[Paving Problem]
Fix $\eps \in (0, 1)$.  Is there a constant $m = m(\eps)$ so that, for sufficiently large $n$, every hollow $n \times n$ matrix has an $(m, \eps)$-paving?
\end{problem}

Anderson \cite{And79:Extreme-Points} has shown that the Paving Problem is equivalent with the Kadison--Singer problem, a major open question in operator theory.  It is closely related to significant problems in harmonic analysis and other areas of mathematics and engineering.  See \cite{CT06:Kadison-Singer-Problem} for an intriguing discussion.

At present, the strongest results on the paving problem are due to Bourgain and Tzafiri \cite{BT91:Problem-Kadison-Singer}.  For a fixed $\eps$, they established that

\vspace{0.5pc}

\begin{enumerate} \setlength{\itemsep}{0.5pc}
\item	every hollow matrix of size $n \times n$ can be paved with at most $m = \bigO(\log n)$ blocks and

\item	every square matrix whose entries are relatively small compared with its dimension can be paved with a constant number of blocks.
\end{enumerate}

\vspace{0.5pc}

\noindent
Let us present a precise statement of their second result.  We use the notation $\oneton{n} \defby \{1, 2, \dots, n\}$.

\begin{thm}[Bourgain--Tzafriri] \label{thm:bdd-paving}
Fix $\gamma > 0$ and $\eps \in (0, 1)$.  There exists a positive integer $m = m(\gamma, \eps)$ so that, for all $n \geq N(\gamma, \eps)$, the following statement holds.  Suppose that $\mtx{A}$ is an $n \times n$ unit-norm matrix with uniformly bounded entries:
$$
\abs{ a_{jk} } \leq \frac{1}{(\log n)^{1 + \gamma}}
\qquad\text{for $j, k = 1, 2, \dots, n$}.
$$
Then there is a partition of the set $\oneton{n}$ into $m$ blocks $\{ \sigma_1, \sigma_2, \dots, \sigma_m \}$ such that
$$
\norm{ \sum\nolimits_{j=1}^m
	\mtx{P}_{\sigma_j} \mtx{A} \mtx{P}_{\sigma_j} }
	\leq \eps
$$
where $\mtx{P}_{\sigma_j}$ is the projector onto the coordinates listed in $\sigma_j$.  The number $m$ satisfies the bound
$$
m \leq \cnst{C} \eps^{- \cnst{C} / \min\{1, \gamma\} }
$$
where $\cnst{C}$ is a positive universal constant.
\end{thm}

The proof of this result published in \cite{BT91:Problem-Kadison-Singer} hinges on a long and delicate calculation of the supremum of a random process.  This computation involves a difficult metric entropy bound based on some subtle iteration arguments.

This note shows that the central step in the known proof can be replaced by another approach based on symmetrization and noncommutative Khintchine inequalities.  This method for studying random matrices is adapted from Rudelson's article \cite{Rud99:Random-Vectors}.  Even though it is simple and elegant, it leads to sharp bounds in many cases.  By itself, Rudelson's technique is not strong enough, so we must also also invoke a method from Bourgain and Tzafrari's proof to complete the argument.  As we go along, we indicate the provenance of various parts of the argument.

\section{Problem Simplifications}

Let us begin with some problem simplifications.  The reductions in this section were all proposed by Bourgain and Tzafriri; we provide proofs for completeness.

The overall strategy is to construct the paving with probabilistic tools.  The first proposition shows that we can leverage a moment estimate for the norm of a random submatrix to build a paving.  The idea is to permute the coordinates randomly and divide them into blocks.  The moment bound shows that, if we restrict the matrix to the coordinates in a random block, then it has small spectral norm.

\begin{prop}[Random Paving Principle]
Fix an integer $m$, and let $n = km$ for an integer $k$.  Let $\mtx{A}$ be an $n \times n$ unit-norm matrix, and suppose that $\mtx{P}$ is a projector onto exactly $k$ coordinates, chosen uniformly at random from the set $\oneton{n}$.  If, for $p \geq \log n$, we have the estimate
$$
\left( \Expect \norm{ \mtx{P} \mtx{A} \mtx{P} }^p \right)^{1/p}
	\leq \eps,
$$
then there exists a partition of the set $\oneton{n}$ into $m$ blocks $\{\sigma_1, \sigma_2, \dots, \sigma_m\}$, each of size $k$, such that
$$
\norm{ \sum\nolimits_{j=1}^m
	\mtx{P}_{\sigma_j} \mtx{A} \mtx{P}_{\sigma_j} }
	\leq 3 \eps
$$
where $\mtx{P}_{\sigma_j}$ is the projector onto the coordinates listed in $\sigma_j$.
\end{prop}


\begin{proof}
Consider a random permutation $\pi$ of the set $\oneton{n}$.  For $j = 1, 2, \dots, m$, define
$$
\sigma_j(\pi) = \{ \pi(jk - k + 1), \pi( jk - k + 2 ), \dots, \pi(jk) \}.
$$
For each $j$, the projector $\mtx{P}_{\sigma_j(\pi)}$ onto the coordinates in $\sigma_j(\pi)$ is a restriction to $k$ coordinates, chosen uniformly at random.  The hypothesis implies that
$$
\Expect \max\nolimits_{j = 1, 2, \dots, m} \smnorm{}{ \mtx{P}_{\sigma_j(\pi)} \mtx{A} \mtx{P}_{\sigma_j(\pi)} }^p \leq m \eps^p.
$$
There must exist a permutation $\pi_0$ for which the left-hand side is smaller than its expectation.  For the partition with blocks $\sigma_j = \sigma_j(\pi_0)$, we have
$$
\norm{ \sum\nolimits_{j = 1}^m \mtx{P}_{\sigma_j} \mtx{A} \mtx{P}_{\sigma_j} }
	=
\max\nolimits_j \norm{ \mtx{P}_{\sigma_j} \mtx{A} \mtx{P}_{\sigma_j} }
	\leq m^{1/p} \eps.
$$
The equality holds because the coordinate blocks are disjoint.  Finally, we have $m^{1/p} \leq \econst$ because $m \leq n$ and $p \geq \log n$.
\end{proof}

This proposition shows that it is sufficient to study the restriction to a random set of coordinates of size $k$.  Although this dependent coordinate model is conceptually simple, it would lead to severe inconveniences later in the proof.  We prefer instead to study an independent coordinate model for the projector where the \emph{expected number} of coordinates equals $k$.  Fortunately, the two models are equivalent for our purposes.

\begin{prop}[Random Coordinate Models] \label{prop:model-equiv}
Fix an integer $m$, and let $n = km$ for an integer $k$.  Let $\mtx{A}$ be an $n \times n$ matrix.  Suppose that $\mtx{P}$ is a projector onto $k$ coordinates, chosen uniformly at random from $\oneton{n}$, and suppose that $\mtx{R}$ is a projector onto a random set of coordinates from $\oneton{n}$, where each coordinate appears independently with probability $k/n$.  For $p > 0$, it holds that
$$
\left( \Expect \norm{ \mtx{P} \mtx{A} \mtx{P} }^p \right)^{1/p}
\leq \left( 2 \Expect \norm{ \mtx{R} \mtx{A} \mtx{R} }^p \right)^{1/p}.
$$
\end{prop}

The reduction to the independent coordinate model also appears in Bourgain and Tzafriri's paper with a different proof.  The following attractive argument is drawn from \cite[Sec.\ 3]{CR06:Quantitative-Robust}.

\begin{proof}
For a coordinate projector $\mtx{R}$, denote by $\sigma(\mtx{R})$ the set of coordinates onto which it projects.  We can make the following computation:
\begin{align*}
\Prob{ \norm{ \mtx{R} \mtx{A} \mtx{R} }^p > t }
	&\geq \sum\nolimits_{j = k}^n \Prob{ \norm{ \mtx{R} \mtx{A} \mtx{R} }^p > t
		\ | \ \# \sigma( \mtx{R} ) = j }
		\cdot \Prob{ \# \sigma(\mtx{R}) = j } \\
	&\geq \Prob{ \norm{ \mtx{R} \mtx{A} \mtx{R} }^p > t
		\ | \ \# \sigma( \mtx{R} ) = k } \cdot
		\sum\nolimits_{j = k}^n \Prob{ \# \sigma(\mtx{R}) = j } \\
	&\geq \frac{1}{2} \Prob{ \norm{ \mtx{P} \mtx{A} \mtx{P} }^p >t }.
\end{align*}
The second inequality holds because the spectral norm of a submatrix is smaller than the spectral norm of the matrix.  The third inequality relies on the fact \cite[Thm.\ 3.2]{JS68:Monotone-Convergence} that the medians of the binomial distribution $\textsc{binomial}( k/n, n )$ lie between $k - 1$ and $k$.  Integrate with respect to $t$ to complete the argument.
\end{proof}

\section{The Main Result}

On account of these simplifications, it suffices to prove the following theorem.  In the sequel, $\mtx{R}_{\delta}$ denotes a square, diagonal matrix whose diagonal entries are independent and identically distributed 0--1 random variables with common mean $\delta$.  The dimensions of $\mtx{R}_\delta$ conform to its context.

\begin{thm} \label{thm:main}
Fix $\gamma > 0$ and $\eps \in (0, 1)$.  There exists a positive integer $m = m(\gamma, \eps)$ so that, for all $n \geq N(\gamma, \eps)$, the following statement holds.  Suppose that $\mtx{A}$ is an $n \times n$ unit-norm matrix with uniformly bounded entries:
$$
\abs{ a_{jk} } \leq \frac{1}{(\log n)^{1 + \gamma}}
\qquad\text{for $j, k = 1, 2, \dots, n$}.
$$
Let $\delta = 1 / m$.  For $p = 2 \cdot \lceil \log n \rceil$, we have
\begin{equation} \label{eqn:main-moment}
\left(\Expect \norm{ \mtx{R}_{\delta} \mtx{A} \mtx{R}_{\delta} }^p \right)^{1/p}
	\leq \eps.
\end{equation}
The number $m$ satisfies the bound
$$
m \leq (0.01 \eps)^{- 2 (1 + \gamma) / \gamma }.
$$
\end{thm}

An example of Bourgain and Tzafriri shows that the number $\gamma$ cannot be removed from the bound $(\log n)^{-(1+\gamma)}$ on the matrix entries \cite[Ex.\ 2.2]{BT91:Problem-Kadison-Singer}.  Fix $\delta \in (0, 1)$.  For each $n \geq N(\delta)$, they exhibit an $n \times n$ matrix $\mtx{A}$ with unit norm and bounded entries:
$$
\abs{ a_{jk} } \leq \frac{2 \log(1 / \delta)}{\log n}.
$$
For this matrix, $\Expect \norm{ \mtx{R}_\delta \mtx{A} \mtx{R}_\delta } \geq 1/2$.  In particular, it has no constant-size random paving when $\eps$ is small.

\begin{proof}[Proof of Theorem \ref{thm:bdd-paving} from Theorem \ref{thm:main}]
Fix $\gamma$ and $\eps$.  Let $m$ be the integer guaranteed by Theorem \ref{thm:main}, and assume that $n$ is sufficiently large.  Suppose we are given an $n \times n$ matrix with unit norm and uniformly bounded entries.  If necessary, augment the matrix with zero rows and columns so that its dimension is a multiple of $m$.

Apply Proposition \ref{prop:model-equiv} to transfer the estimate \eqref{eqn:main-moment} to the dependent coordinate model.  The Random Paving Principle shows that the augmented matrix has an $(m, 6\eps)$-paving.  Discard the zero rows and columns to complete the proof of Theorem \ref{thm:bdd-paving}.
\end{proof}

\section{Proof of Theorem \ref{thm:main}}

In this section, we establish Theorem \ref{thm:main}.  The proofs of the supporting results are postponed to the subsequent sections.

Fix $\gamma > 0$ and $\eps \in (0, 1)$.  We assume for convenience that $n \geq 8$, and we suppose that $\mtx{A}$ is an $n \times n$ matrix with unit norm and uniformly bounded entries:
$$
\abs{ a_{jk} } \leq \frac{1}{(\log n)^{1 + \gamma}} \defby \mu.
$$
In the sequel, the symbol $\mu$ always abbreviates the uniform bound.  Finally, set $p = 2 \cdot \lceil \log n \rceil$.

The major task in the proof is to obtain an estimate for the quantity
$$
E(\varrho) \defby \left( \Expect \norm{ \mtx{R}_\varrho \mtx{A}\mtx{R}_\varrho }^p \right)^{1/p},
$$
where $\varrho$ is not too small.  This estimate is accomplished with decoupling, symmetrization, and noncommutative Khintchine inequalities.  This approach is adapted from work of Rudelson \cite{Rud99:Random-Vectors} and Rudelson--Vershynin \cite{RV07:Sampling-Large}.  Given this estimate for $E(\varrho)$, we extrapolate the value of $E( m^{-1} )$ for a large constant $m = m(\gamma, \eps)$.  This step relies on an elegant method due to Bourgain and Tzafriri.

Before continuing, we instate a few more pieces of notation.  The symbol $\pnorm{1, 2}{ \cdot }$ denotes the norm of an operator mapping $\ell_1^n$ to $\ell_2^n$.  For a matrix $\mtx{X}$ expressed in the standard basis, $\pnorm{1, 2}{ \mtx{X} }$ is the maximum $\ell_2^n$ norm achieved by a column of $\mtx{X}$.  The norm $\pnorm{\max}{\mtx{X}}$ calculates the maximum absolute value of an entry of $\mtx{X}$.

\subsection{Step 1: Decoupling}

As in Bourgain and Tzafriri's work, the first step involves a classical decoupling argument.  First, we must remove the diagonal of the matrix.  Since the entries of $\mtx{A}$ do not exceed $\mu$, it follows that $\norm{ \diag \mtx{A} } \leq \mu$.  Define
$$
\mtx{B} = \frac{1}{1 + \mu} (\mtx{A} - \diag{\mtx{A}}).
$$
Note that $\mtx{B}$ has a zero diagonal and that $\norm{\mtx{B}} \leq 1$.  Furthermore,
$$
\abs{b_{jk}} < \mu \qquad\text{for $j, k = 1, 2, \dots, n$}.
$$
With this definition,
$$
E( \varrho ) \leq \norm{ \diag{\mtx{A}} } + (1 + \mu) \left( \Expect
	\norm{ \mtx{R}_\varrho \mtx{B}
		\mtx{R}_\varrho }^p \right)^{1/p}.
$$
The expectation on the right-hand side cannot exceed one, so we have
$$
E( \varrho ) \leq 2 \mu + \left( \Expect
	\norm{ \mtx{R}_\varrho \mtx{B} \mtx{R}_\varrho }^p \right)^{1/p}.
$$
Now, we may replace the projector $\mtx{R}_\varrho$ by a pair of independent projectors by invoking the following result.

\begin{prop} \label{prop:decoupling}
Let $\mtx{B}$ be a square matrix with a zero diagonal, and let $p \geq 1$.  Then
$$
\left( \Expect \norm{ \mtx{R}_\varrho \mtx{B} \mtx{R}_\varrho }^p \right)^{1/p}
	\leq 20 \left( \Expect \smnorm{}{ \mtx{R}_{\varrho} \mtx{B}
	\mtx{R}_{\varrho}' }^p \right)^{1/p}
$$
where the two random projectors on the right-hand side are independent.
\end{prop}

\noindent
See \cite[Prop.\ 1.1]{BT87:Invertibility-Large} or \cite[Sec.\ 4.4]{LT91:Probability-Banach} for the simple proof.

We apply Proposition \ref{prop:decoupling} to reach
\begin{equation} \label{eqn:step1-end}
E( \varrho ) \leq 2 \mu + 20 \left(
	\Expect \smnorm{}{ \mtx{R}_{\varrho} \mtx{B}
		\mtx{R}_{\varrho}' }^p \right)^{1/p}.
\end{equation}

\subsection{Step 2: Norm of a Random Restriction}

The next step of the proof is to develop a bound on the spectral norm of a matrix that has been restricted to a random subset of its columns.  The following result is due to Rudelson and Vershynin \cite{RV07:Sampling-Large}, with some inessential modifications by the author.

\begin{prop}[Rudelson--Vershynin] \label{prop:random-restrict}
Let $\mtx{X}$ be an $n \times n$ matrix, and suppose that $p \geq 2 \log n \geq 2$.  Then
$$
\left( \Expect \norm{\mtx{X} \mtx{R}_\varrho}^p \right)^{1/p}
	\leq 3 \sqrt{p} \left(
		\Expect \pnorm{1,2}{ \mtx{X}\mtx{R}_\varrho }^p \right)^{1/p}
		+ \sqrt{\varrho} \norm{\mtx{X}}.
$$
\end{prop}

The proof of Proposition \ref{prop:random-restrict} depends on a lemma of Rudelson that bounds the norm of a Rademacher sum of rank-one, self-adjoint matrices \cite{Rud99:Random-Vectors}.  This lemma, in turn, hinges on the noncommutative Khintchine inequality \cite{L-P86:Inegalites-Khintchine,Buc01:Operator-Khintchine}.  See Section \ref{app:restriction} for the details.

To account for the influence of $\mtx{R}_\varrho'$, we apply Proposition \ref{prop:random-restrict} with $\mtx{X} = \mtx{R}_\varrho \mtx{B}$.  Inequality \eqref{eqn:step1-end} becomes
$$
E(\varrho) \leq 2 \mu + 60 \sqrt{p} \left( \Expect
	\pnorm{1,2}{ \mtx{R}_\varrho \mtx{B} \mtx{R}_\varrho' }^p \right)^{1/p}
		+ 20 \sqrt{\varrho} \left( \Expect
		\norm{\mtx{R}_\varrho \mtx{B}}^p \right)^{1/p}.
$$
We invoke Proposition \ref{prop:random-restrict} again with $\mtx{X} = \mtx{B}^\adj$ to reach
$$
E(\varrho) \leq 2 \mu + 60 \sqrt{p} \left( \Expect
	\pnorm{1,2}{ \mtx{R}_\varrho \mtx{B} \mtx{R}_\varrho' }^p \right)^{1/p}
	+ 60 \sqrt{\varrho p} \left( \Expect
	\pnorm{1,2}{ \mtx{B}^\adj \mtx{R}_\varrho }^p \right)^{1/p}
	+ 20 \varrho \norm{\mtx{B}^\adj}.
$$
Discard the projector $\mtx{R}_\varrho'$ from the first expectation by means of the observation
$$
\pnorm{1,2}{ \mtx{R}_\varrho \mtx{B} \mtx{R}_\varrho' }
	\leq \pnorm{1,2}{\mtx{R}_\varrho \mtx{B}}.
$$
In words, the maximum column norm of a matrix exceeds the maximum column norm of any submatrix.  We also have the bound
$$
\pnorm{1,2}{\mtx{B}^\adj \mtx{R}_\varrho }
	\leq \pnorm{1, 2}{ \mtx{B}^\adj }
	\leq \norm{ \mtx{B}^\adj }
	\leq 1
$$
because the spectral norm dominates the maximum $\ell_2^n$ norm of a column.  The inequality $\varrho \leq \sqrt{\varrho}$ yields
\begin{equation} \label{eqn:step2-end}
E(\varrho) \leq 2\mu + 60 \sqrt{p} \left( \Expect
	\pnorm{1,2}{ \mtx{R}_\varrho \mtx{B} }^p \right)^{1/p}
	+ 80 \sqrt{\varrho p}.
\end{equation}

\subsection{Step 3: Estimate of Maximum Column Norm}

To complete our estimate of $E(\varrho)$, we must bound the remaining expectation.  The following result does the job.

\begin{prop} \label{prop:12-norm}
Let $\mtx{X}$ be an $n \times n$ matrix, and suppose that $p \geq 2 \log n \geq 4$.  Then
$$
\left( \Expect \pnorm{1,2}{\mtx{R}_\varrho \mtx{X}}^p \right)^{1/p}
	\leq 3 \sqrt{p} \pnorm{\max}{\mtx{X}}
		+ \sqrt{\varrho} \pnorm{1,2}{\mtx{X}}.
$$
\end{prop}

The proof of Proposition \ref{prop:12-norm} uses only classical methods, including symmetrization and scalar Khintchine inequalities.  A related bound appears inside the proof of Proposition 2.5 in \cite{BT91:Problem-Kadison-Singer}.  Turn to Section \ref{app:12-norm} for the argument.

Apply Proposition \ref{prop:12-norm} to the remaining expectation in \eqref{eqn:step2-end} to find that
$$
E(\varrho) \leq 2 \mu + 180 p \pnorm{\max}{\mtx{B}} + 60
	\sqrt{\varrho p} \pnorm{1,2}{\mtx{B}} + 80 \sqrt{\varrho p}.
$$
As above, the maximum column norm $\pnorm{1,2}{ \mtx{B} } \leq 1$.
The entries of $\mtx{B}$ are uniformly bounded by $\mu$.  Recall $p = 2 \cdot \lceil \log n \rceil$ to conclude that
\begin{equation} \label{eqn:step3-end}
E(\varrho) \leq 550 \mu \log n + 250 \sqrt{\varrho \log n},
\end{equation}
taking into account $\lceil \log n \rceil \leq 1.5 \log n$ whenever $n \geq 8$.

The result in \eqref{eqn:step3-end} is not quite strong enough to establish Theorem \ref{thm:main}.  In the theorem, the relation between the size $m$ of the paving and the proportion $\delta$ of columns is $\delta = 1/m$.  The parameter $\varrho$ also represents the proportion of columns selected.  Unfortunately, when we set $\varrho = 1/m$, we find that the bound in \eqref{eqn:step3-end} is trivial unless $\varrho$ is smaller than $\cnst{c} / \log n$, which suggests that $m$ grows logarithmically with $n$.  To prove the result, however, we must obtain a bound for $m$ that is independent of dimension.

\subsection{Step 4: Extrapolation}

To finish the argument, we require a remarkable fact uncovered by Bourgain and Tzafriri in their work.  Roughly speaking, the value of $E(\varrho)^p$ is comparable with a polynomial of low degree.  It is possible to use the inequality \eqref{eqn:step3-end} to estimate the coefficients of this polynomial.  We can then extrapolate to obtain a nontrivial estimate of $E(\delta)^p$, where $\delta$ is a small constant.

\begin{prop}[Bourgain--Tzafriri] \label{prop:extrap}
Let $\mtx{X}$ be an $n \times n$ matrix with $\norm{\mtx{X}} \leq 1$.  Suppose that $p$ is an even integer with $p \geq 2 \log n$.  Choose parameters $\delta \in (0, 1)$ and $\varrho \in (0, 0.5)$.  For each $\lambda \in (0, 1)$, it holds that
$$
\left( \Expect \norm{ \mtx{R}_\delta \mtx{X}\mtx{R}_\delta }^p \right)^{1/p}
	\leq 60 \left[ \delta^{\lambda} + \varrho^{-\lambda}
	\left( \Expect \norm{ \mtx{R}_{\varrho} \mtx{X}\mtx{R}_{\varrho} }^p \right)^{1/p} \right].
$$
\end{prop}

The proof depends essentially on a result of V.\ A.\ Markov that bounds the coefficients of a polynomial in terms of its maximum value.  See Section \ref{sec:extrap} for the details.

Recall now that
$$
\mu \leq \frac{1}{(\log n)^{1 + \gamma}}.
$$
Set the proportion $\varrho = (\log n)^{-1 - 2\gamma}$, and introduce these quantities into \eqref{eqn:step3-end} to obtain
$$
\left( \Expect
	\norm{ \mtx{R}_\varrho \mtx{A}\mtx{R}_\varrho }^p \right)^{1/p}
	\leq 800 (\log n)^{- \gamma}.
$$
Proposition \ref{prop:extrap} shows that
$$
\left( \Expect \norm{ \mtx{R}_\delta \mtx{A}\mtx{R}_\delta }^p \right)^{1/p}
	\leq 60 \delta^{\lambda} + 48000 (\log n)^{\lambda (1 + 2\gamma) - \gamma}
$$
for every value of $\lambda$ in $(0, 1)$.  Make the selection $\lambda = \gamma / (2 + 2\gamma)$.  Since the exponent on the logarithm is strictly negative, it follows for sufficiently large $n$ that
$$
\left( \Expect \norm{ \mtx{R}_\delta \mtx{A}\mtx{R}_\delta }^p \right)^{1/p} \leq 100 \, \delta^{\gamma / (2 + 2\gamma)}.
$$
To make the right-hand side less than a parameter $\eps$, it suffices that $\delta \leq (0.01 \eps)^{2 (1 + \gamma)/\gamma}$.  Therefore, any value
$$
m \geq (0.01 \eps)^{- 2(1 + \gamma) / \gamma}
$$
is enough to establish Theorem \ref{thm:main}.




\section{Proof of Random Restriction Estimate} \label{app:restriction}

In this section, we establish Proposition \ref{prop:random-restrict}.  The difficult part of the estimation is performed with the noncommutative Khintchine inequality.  This result was originally discovered by Lust-Picquard \cite{L-P86:Inegalites-Khintchine}. We require a sharp version due to Buchholz \cite{Buc01:Operator-Khintchine} that provides the optimal order of growth in the constant.

Before continuing, we state a few definitions.  Given a matrix $\mtx{X}$, let $\sigma( \mtx{X} )$ denote the vector of its singular values, listed in weakly decreasing order.  The \term{Schatten $p$-norm} $\pnorm{S_p}{\cdot}$ is calculated as
$$
\pnorm{S_p}{\mtx{X}} = \pnorm{p}{\vct{\sigma}}
$$
where $\pnorm{p}{\cdot}$ denotes the $\ell_p$ vector norm.

A \term{Rademacher random variable} takes the two values $\pm 1$ with equal probability.  A \term{Rademacher sequence} is a sequence of independent Rademacher variables.

\begin{prop}[Noncommutative Khintchine Inequality] \label{prop:khintchine}
Let $\{ \mtx{X}_j \}$ be a finite sequence of matrices of the same dimension, and let $\{\eps_j\}$ be a Rademacher sequence.  For each $p \geq 2$,
\begin{equation} \label{eqn:khintchine}
\left[ \Expect
	\pnorm{S_{p}}{ \sum\nolimits_j \eps_j \mtx{X}_j }^{p} \right]^{1/p}
	\leq \cnst{C}_{p} \max\left\{
	\pnorm{S_{p}}{ \left( \sum\nolimits_j \mtx{X}_j \mtx{X}_j^\adj \right)^{1/2} },
	\pnorm{S_{p}}{ \left( \sum\nolimits_j \mtx{X}_j^\adj \mtx{X}_j \right)^{1/2} } \right\},
\end{equation}
where $\cnst{C}_{p} \leq 2^{-0.25} \sqrt{\pi/\econst}\, \sqrt{p}$.
\end{prop}

This proposition is a corollary of Theorem 5 of \cite{Buc01:Operator-Khintchine}.  In this work, Buchholz shows that the noncommutative Khintchine inequality holds with a Gaussian sequence in place of the Rademacher sequence.  He computes the optimal constant when $p$ is an even integer:
$$
\cnst{C}_{2n} = \left( \frac{(2n)!}{2^n n!} \right)^{1/2n}.
$$
One extends this result to other values of $p$ using Stirling's approximation and an interpolation argument.  The inequality for Rademacher variables follows from the contraction principle.

In an important paper \cite{Rud99:Random-Vectors}, Rudelson showed how to use the noncommutative Khintchine inequality to study the moments of a Rademacher sum of rank-one matrices.

\begin{lemma}[Rudelson] \label{lem:rudelson}
Suppose that $\vct{x}_1, \vct{x}_2, \dots, \vct{x}_n$ are the columns of a matrix $\mtx{X}$.  For any $p \geq 2 \log n$, it holds that
$$
\left( \Expect \norm{ \sum\nolimits_{j=1}^n \eps_j \vct{x}_j \vct{x}_j^\adj }^p \right)^{1/p}
	\leq 1.5 \sqrt{p}
	\pnorm{1,2}{ \mtx{X} } \norm{ \mtx{X} },
$$
where $\{\eps_j\}$ is a Rademacher sequence.
\end{lemma}

\begin{proof}
First, bound the spectral norm by the Schatten $p$-norm.
$$
E \defby
\left( \Expect \norm{ \sum\nolimits_{j=1}^n \eps_j \vct{x}_j \vct{x}_j^\adj }^p \right)^{1/p}
	\leq \left( \Expect \pnorm{S_p}{ \sum\nolimits_{j=1}^n \eps_j \vct{x}_j \vct{x}_j^\adj }^p \right)^{1/p}.
$$
Apply the noncommutative Khintchine inequality to obtain
$$
E \leq
\cnst{C}_p \pnorm{S_p}{ \left( \sum\nolimits_{j=1}^n
	\enormsq{\vct{x}_j} \vct{x}_j \vct{x}_j^\adj \right)^{1/2}}.
$$
The rank of matrix inside the norm does not exceed $n$, so we can bound the Schatten $p$-norm by the spectral norm if we pay a factor of $n^{1/p}$, which does not exceed $\sqrt{\econst}$.  Afterward, pull the square root out of the norm to find
$$
E \leq \cnst{C}_p \sqrt{\econst} \norm{ \sum\nolimits_{j=1}^n
	\enormsq{\vct{x}_j} \vct{x}_j \vct{x}_j^\adj }^{1/2}.
$$
The summands are positive semidefinite, so the spectral norm of the sum increases monotonically with each scalar coefficient.  Therefore, we may replace each coefficient by $\max_j \enormsq{\vct{x}_j}$ and use the homogeneity of the norm to obtain
$$
E \leq
\cnst{C}_p \sqrt{\econst} \max\nolimits_j \enorm{\vct{x}_j}
	\norm{ \sum\nolimits_{j=1}^n \vct{x}_j \vct{x}_j^\adj }^{1/2}.
$$
The maximum can be rewritten as $\pnorm{1,2}{\mtx{X}}$, and the spectral norm can be expressed as
$$
\norm{ \sum\nolimits_{j=1}^n \vct{x}_j \vct{x}_j^\adj }^{1/2}
	= \norm{ \mtx{X}\mtx{X}^\adj }^{1/2} = \norm{ \mtx{X} }.
$$
Recall that $\cnst{C}_p \leq 2^{-0.25} \sqrt{\pi/\econst}\, \sqrt{p}$ to complete the proof.
\end{proof}

Recently, Rudelson and Vershynin showed how Lemma \ref{lem:rudelson} implies a bound on the moments of the norm of a matrix that is compressed to a random subset of columns \cite{RV07:Sampling-Large}.

\begin{prop}[Rudelson--Vershynin] \label{prop:random-restriction}
Let $\mtx{X}$ be a matrix with $n$ columns, and suppose that $p \geq 2 \log n \geq 2$.  It holds that
$$
\left( \Expect \norm{ \mtx{X} \mtx{R}_{\varrho} }^{p} \right)^{1/p}
	\leq 3 \sqrt{p}
	\left( \Expect \pnorm{1,2}{ \mtx{X} \mtx{R}_\varrho }^{p} \right)^{1/p}
		+ \sqrt{\varrho} \norm{ \mtx{X} }.
$$
\end{prop}

In words, a random compression of a matrix gets its share of the spectral norm plus another component that depends on the total number of columns and on the $\ell_2^n$ norms of the columns.

\begin{proof}
Let us begin with an overview of the proof.  First, we express the random compression as a random sum.  Then we symmetrize the sum and apply Rudelson's lemma to obtain an upper bound involving the value we are trying to estimate.  Finally, we solve an algebraic relation to obtain an explicit estimate for the moment.

We seek a bound for
$$
E \defby \left( \Expect \norm{ \mtx{X} \mtx{R}_\varrho }^{p} \right)^{1/p}.
$$
First, observe that
$$
E^2 =
\left( \Expect \norm{ \mtx{X} \mtx{R}_\varrho \mtx{X}^\adj }^{p/2} \right)^{2/p}
= \left( \Expect \norm{ \sum\nolimits_{j=1}^n \varrho_j \vct{x}_j \vct{x}_j^\adj }^{p/2} \right)^{2/p}
$$
where $\{ \varrho_j \}$ is a sequence of independent 0--1 random variables with common mean $\varrho$.  Subtract the mean, and apply the triangle inequality (once for the spectral norm and once for the $L_{p/2}$ norm):
$$
E^2 \leq \left( \Expect \norm{ \sum\nolimits_{j=1}^n (\varrho_j - \varrho) \vct{x}_j \vct{x}_j^\adj }^{p/2} \right)^{2/p}
	+ \varrho \norm{ \sum\nolimits_{j=1}^n \vct{x}_j \vct{x}_j^\adj }.
$$
In the sum, write $\varrho = \Expect \varrho_j'$ where $\{ \varrho_j' \}$ is an independent copy of the sequence $\{\varrho_j \}$.  Draw the expectation out of the norm with Jensen's inequality:
$$
E^2 \leq \left( \Expect \norm{ \sum\nolimits_{j=1}^n (\varrho_j - \varrho_j') \vct{x}_j \vct{x}_j^\adj }^{p/2} \right)^{2/p}
	+ \varrho \norm{ \mtx{X}\mtx{X}^\adj }.
$$
The random variables $(\varrho_j - \varrho_j')$ are symmetric and independent, so we may symmetrize them using the standard method, Lemma 6.1 of \cite{LT91:Probability-Banach}.
$$
E^2 \leq \left( \Expect \norm{ \sum\nolimits_{j=1}^n \eps_j (\varrho_j - \varrho_j') \vct{x}_j \vct{x}_j^\adj }^{p/2} \right)^{2/p}
	+ \varrho \norm{ \mtx{X} }^2
$$
where $\{\eps_j\}$ is a Rademacher sequence.  Apply the triangle inequality again and use the identical distribution of the sequences to obtain
$$
E^2 \leq 2 \left( \Expect \norm{ \sum\nolimits_{j=1}^n \eps_j \varrho_j \vct{x}_j \vct{x}_j^\adj }^{p/2} \right)^{2/p}
	+ \varrho \norm{ \mtx{X} }^2
$$
Writing $\Omega = \{ j : \varrho_j = 1 \}$, we see that
$$
E^2 \leq 2 \left[ \Expect_\Omega \left( \Expect_{\vct{\eps}} \norm{ \sum\nolimits_{\Omega} \eps_j \vct{x}_j \vct{x}_j^\adj }^{p/2} \right)^{(2/p) (p/2)} \right]^{2/p}
	+ \varrho \norm{ \mtx{X} }^2.
$$
Here, $\Expect_{\vct{\eps}}$ is the partial expectation with respect to $\{ \eps_j \}$, holding the other random variables fixed.

To estimate the large parenthesis, invoke Rudelson's Lemma, conditional on $\Omega$.  The matrix in the statement of the lemma is $\mtx{X} \mtx{R}_\varrho$, resulting in
$$
E^2 \leq 3 \sqrt{p}
	\left[ \Expect \left( \pnorm{1,2}{ \mtx{X} \mtx{R}_\varrho } \norm{ \mtx{X} \mtx{R}_\varrho } \right)^{p/2} \right]^{2/p}
	+ \varrho \norm{ \mtx{X} }^2.
$$
Apply the Cauchy--Schwarz inequality to find that
$$
E^2 \leq 3 \sqrt{p}
	\left( \Expect \pnorm{1,2}{ \mtx{X} \mtx{R}_\varrho }^{p} \right)^{1/p}
	\left( \Expect \norm{ \mtx{X} \mtx{R}_\varrho }^{p} \right)^{1/p}
	+ \varrho \norm{ \mtx{X} }^2.
$$

This inequality takes the form $E^2 \leq b E + c$.  Select the larger root of the quadratic and use the subadditivity of the square root:
$$
E \leq \frac{ b + \sqrt{ b^2 + 4c } }{2}
	\leq b + \sqrt{c}.
$$
This yields the conclusion.
\end{proof}

\section{Proof of Maximum Column Norm Estimate} \label{app:12-norm}

This section establishes the moment bound for the maximum column norm of a matrix that has been restricted to a random set of its rows.  We use an approach that is analogous with the argument in Proposition \ref{prop:random-restriction}.  In this case, we require only the scalar Khintchine inequality to perform the estimation.  Bourgain and Tzafriri's proof of Proposition 2.5 \cite{BT91:Problem-Kadison-Singer} contains a similar bound, developed with a similar argument.

\begin{prop}
Assume that $\mtx{X}$ has $n$ columns, and suppose $p \geq 2 \log n \geq 4$.  Then
$$
\left( \Expect \pnorm{1,2}{ \mtx{R}_\varrho \mtx{X} }^{p} \right)^{1/p}
	\leq 2^{1.5} \sqrt{p} \pnorm{\max}{\mtx{X}}
		+ \sqrt{\varrho} \pnorm{1,2}{\mtx{X}}.
$$
\end{prop}

In words, the $B(\ell_1^n, \ell_2^n)$ norm of a matrix that has been compressed to a random set of rows gets its share of the total, plus an additional component that depends on the number of columns and the magnitude of the largest entry in matrix.

\begin{proof}
Our strategy is the same as in the proof of Proposition \ref{prop:random-restriction}, so we pass lightly over the details.  Let $\{\varrho_j\}$ be a sequence of independent 0--1 random variables with common mean $\varrho$.  We seek a bound for
$$
E^2 \defby \left( \Expect
	\pnorm{1,2}{ \mtx{R}_\varrho \mtx{X} }^{p} \right)^{2/p}
	= \left( \Expect \max\nolimits_{k = 1, 2, \dots, n}
		\abs{ \sum\nolimits_j \varrho_j \, \abssq{ x_{jk} } }^{p/2} \right)^{2/p}.
$$
In the sequel, we abbreviate $q = p/2$ and also $y_{jk} = \abssq{x_{jk}}$.

First, center and symmetrize the selectors.  In the following calculation, $\{ \varrho_j' \}$ is an independent copy of the sequence $\{\varrho_j\}$, and $\{\eps_j\}$ is a Rademacher sequence, independent of everything else.
\begin{align*}
E^2 &\leq \left( \Expect \max\nolimits_k  \abs{ \sum\nolimits_j
		(\varrho_j - \varrho) y_{jk} }^{q} \right)^{1/q}
		+ \varrho \max\nolimits_k \abs{ \sum\nolimits_j y_{jk} } \\
	&\leq \left( \Expect \max\nolimits_k \abs{ \sum\nolimits_j
		(\varrho_j - \varrho_j') y_{jk} }^{q} \right)^{1/q}
		+ \varrho \pnorm{1,2}{ \mtx{X} }^2 \\
	&= \left( \Expect \max\nolimits_k \abs{ \sum\nolimits_j
		\eps_j (\varrho_j - \varrho_j') y_{jk} }^{q} \right)^{1/q}
		+ \varrho \pnorm{1,2}{ \mtx{X} }^2 \\
	&\leq 2 \left( \Expect \max\nolimits_k \abs{ \sum\nolimits_j
		\eps_j \varrho_j y_{jk} }^{q} \right)^{1/q}
		+ \varrho \pnorm{1,2}{ \mtx{X} }^2.
\end{align*}
The first step uses the triangle inequality; the second uses $\varrho = \Expect \varrho_j'$ and Jensen's inequality; the third follows from the standard symmetrization, Lemma 6.1 of \cite{LT91:Probability-Banach}.  In the last step, we invoked the triangle inequality and the identical distribution of the two sequences.

Next, bound the maximum by a sum and introduce conditional expectations:
$$	
E^2 \leq 2 \left( \Expect_{\vct{\varrho}} 
		\sum\nolimits_k \Expect_{\vct{\eps}} \abs{ \sum\nolimits_j
		\eps_j \varrho_j y_{jk} }^{q} \right)^{1/q}
		+ \varrho \pnorm{1,2}{ \mtx{X} }^2.
$$
Here, $\Expect_{\vct{\eps}}$ denotes partial expectation with respect to $\{ \eps_j \}$, holding the other random variables fixed.  Since $q \geq 2$, we may apply the scalar Khintchine inequality to the inner expectation to obtain
$$
E^2 \leq 2 \cnst{C}_q \left( \Expect_{\vct{\varrho}} 
		\sum\nolimits_k \abs{ \sum\nolimits_j
		\varrho_j y_{jk}^2 }^{q/2} \right)^{1/q}
		+ \varrho \pnorm{1,2}{ \mtx{X} }^2,
$$
where the constant $\cnst{C}_q \leq 2^{0.25} \econst^{-1/2} \sqrt{q}$.  The value of the constant follows from work of Haagerup \cite{Hag82:Best-Constants}, combined with Stirling's approximation.

Bound the outer sum, which ranges over $n$ indices, by a maximum:
$$
E^2 \leq 2^{1.25} \econst^{-1/2} n^{1/q} \sqrt{q}
	\left( \Expect_{\vct{\varrho}} 
		\max\nolimits_k \abs{ \sum\nolimits_j
		\varrho_j y_{jk}^2 }^{q/2} \right)^{1/q}
		+ \varrho \pnorm{1,2}{ \mtx{X} }^2.
$$
Since $q \geq \log n$, it holds that $n^{1/q} \leq \econst$, which implies that the leading constant is less than four.  Use H{\"o}lder's inequality to bound the sum, and then apply H{\"o}lder's inequality again to double the exponent:
\begin{align*}
E^2 &< 4 \sqrt{q}
	\biggl( \max_{j,k} y_{jk} \biggr)^{1/2} 
	\left( \Expect_{\vct{\varrho}} 
		\max\nolimits_k \abs{ \sum\nolimits_j
		\varrho_j y_{jk} }^{q/2} \right)^{1/q}
		+ \varrho \pnorm{1,2}{ \mtx{X} }^2 \\
	&\leq 4 \sqrt{q}
	\biggl( \max_{j,k} y_{jk} \biggr)^{1/2} 
		\biggl( \Expect_{\vct{\varrho}}
			\max\nolimits_k \abs{ \sum\nolimits_j
			\varrho_j y_{jk} }^{q} \biggr)^{1/2q}
		+ \varrho \pnorm{1,2}{ \mtx{X} }^2.
\end{align*}
Recall that $q = p/2$ and that $y_{jk} = \abssq{x_{jk}}$.  Observe that we have obtained a copy of $E$ on the right-hand side, so
$$
E^2 \leq 2^{1.5} \sqrt{p} \pnorm{\max}{ \mtx{X} } E
		+ \varrho \pnorm{1,2}{ \mtx{X} }^2.
$$

As in the proof of Proposition \ref{prop:random-restriction}, we take the larger root of the quadratic and invoke the subadditivity of the square root to reach
$$
E \leq 2^{1.5} \sqrt{p} \pnorm{\max}{\mtx{X}}
		+ \sqrt{\varrho} \pnorm{1,2}{\mtx{X}}.
$$
This is the advertised conclusion.
\end{proof}

\section{Proof of Extrapolation Bound} \label{sec:extrap}

This section summarizes the argument of Bourgain and Tzafriri that leads to the extrapolation result.  The key to the proof is an observation due to V.\ A.\ Markov that estimates the coefficients of an arbitrary polynomial in terms of its maximum value \cite[Sec.~2.9]{Tim63:Theory-Approximation}.

\begin{prop}[Markov] \label{prop:markov}
Let $r(t) = \sum_{k = 0}^d c_k t^k$.  The coefficients of the polynomial $r$ satisfy the inequality
$$
\abs{c_k} \leq \frac{d^k}{k!} \max_{\abs{t} \leq 1} \abs{ r(t) }
	\leq \econst^d \max_{\abs{t} \leq 1} \abs{ r(t) }.
$$
for each $k = 0, 1, \dots, d$.
\end{prop}

The proof depends on the minimax property of the Chebyshev polynomial of degree $d$, combined with a careful determination of its coefficients.

\begin{prop}[Bourgain--Tzafriri] 
Let $p$ be an even integer with $p \geq 2 \log n$.  Suppose that $\mtx{X}$ is an $n \times n$ matrix with $\norm{\mtx{X}} \leq 1$.  Choose parameters $\delta \in (0, 1)$ and $\varrho \in (0, 0.5)$.  For each $\lambda \in (0, 1)$, it holds that
$$
\left( \Expect \norm{ \mtx{R}_\delta \mtx{X}\mtx{R}_\delta }^p \right)^{1/p}
	\leq 60 \left[ \delta^{\lambda} + \varrho^{-\lambda}
	\left( \Expect \norm{ \mtx{R}_{\varrho} \mtx{X}\mtx{R}_{\varrho} }^p \right)^{1/p} \right].
$$
For self-adjoint matrices, the constant is halved.
\end{prop}

\begin{proof}
We assume that $\mtx{X}$ is self-adjoint.  For general $\mtx{X}$, apply the final bound to each half of the Cartesian decomposition
$$
\mtx{X} = \frac{\mtx{X} + \mtx{X}^\adj}{2}
	+ \frac{ \iunit(\mtx{X} - \mtx{X}^\adj) }{2 \iunit}.
$$
This yields the same result with constants doubled.

Consider the function
$$
F(s) = \Expect \norm{ \mtx{R}_s \mtx{X} \mtx{R}_s }^p
\qquad\text{with $0 \leq s \leq 1$.}
$$
Note that $F(s) \leq 1$ because $\norm{ \mtx{R}_s \mtx{X} \mtx{R}_s } \leq \norm{\mtx{X}} \leq 1$.  Furthermore, $F$ increases monotonically.

Next, we show that $F$ is comparable with a polynomial.  Use the facts that $p$ is even, that $p \geq \log n$, and that $\rank{\mtx{X}} \leq n$ to check the inequalities
$$
F(s) \leq
	\Expect \trace (\mtx{R}_s \mtx{X} \mtx{R}_s)^p
	\leq \econst^p F( s ).
$$
It is easy to see that the central member is a polynomial of maximum degree $p$ in the variable $s$.  Indeed, one may expand the product and compute the expectation using the fact that the diagonal entries of $\mtx{R}_{s}$ are independent 0--1 random variables of mean $s$.  Therefore,
$$
\Expect \trace (\mtx{R}_s \mtx{X} \mtx{R}_s)^p
	= \sum_{k = 1}^p c_k s^k
$$
for (unknown) coefficients $c_1, c_2, \dots, c_p$.  The polynomial has no constant term because $\mtx{R}_0 = \mtx{0}$.

We must develop some information about this polynomial.  Make the change of variables $s = \varrho t^2$ to see that
$$
\abs{ \sum\nolimits_{k=1}^p c_k \varrho^k t^{2k} }
	\leq \econst^p F( \varrho t^2 )
	\leq \econst^p F( \varrho )
\qquad\text{when $\abs{t} \leq 1$.}
$$
The second inequality follows from the monotonicity of $F$.
The polynomial on the left-hand side has degree $2p$ in the variable $t$, so Proposition \ref{prop:markov} results in
$$
\abs{c_k} \varrho^k \leq \econst^{3p} F( \varrho ).
\qquad\text{for $k = 1, 2, \dots, p$.} 
$$
From here, it also follows that $\abs{c_k} \leq \econst^{3p}$ by taking $\varrho = 1$.

Finally, we directly evaluate the polynomial at $\delta$ using the facts we have uncovered.  For an arbitrary value of $\lambda$ in $(0,1)$, we have
\begin{align*}
F(\delta) &\leq
\abs{ \sum\nolimits_{k = 1}^d c_k \delta^k } \\
	&\leq \sum\nolimits_{k=1}^{\lfloor \lambda p \rfloor} \abs{c_k}
		+ \sum\nolimits_{k=1 + \lfloor \lambda p \rfloor}^p
			\abs{c_k} \delta^k \\
	&\leq \econst^{3p} F( \varrho )
		\sum_{k=1}^{\lfloor \lambda p \rfloor} \varrho^{-k}
		+ p \econst^{3p} \delta^{\lambda p} \\
	&\leq 2 \econst^{3p} \varrho^{- \lambda p} F(\varrho)
		+ p \econst^{3p} \delta^{\lambda p} 
\end{align*}
since $\varrho \leq 0.5$.  Since $x \mapsto x^{1/p}$ is subadditive, we conclude that
$$
F( \delta )^{1/p} \leq
	2^{1/p} \econst^3 \cdot \varrho^{-\lambda} F(\varrho)^{1/p}
	+ p^{1/p} \econst^3 \cdot \delta^\lambda.
$$
A numerical calculation shows that both the leading terms are less than 30, irrespective of $p$.
\end{proof}

\section*{Acknowledgments}

I wish to thank Roman Vershynin for encouraging me to study the paving problem.

\bibliographystyle{alpha}
\bibliography{paving}

\begin{thebibliography}{And79}

\bibitem[And79]{And79:Extreme-Points}
J.~Anderson.
\newblock Extreme points in sets of positive linear maps on {$B(H)$}.
\newblock {\em J. Functional Anal.}, 31:195--217, 1979.

\bibitem[BT87]{BT87:Invertibility-Large}
J.~Bourgain and L.~Tzafriri.
\newblock Invertibility of ``large'' submatrices with applications to the
  geometry of {B}anach spaces and harmonic analysis.
\newblock {\em Israel J. Math}, 57(2):137--224, 1987.

\bibitem[BT91]{BT91:Problem-Kadison-Singer}
J.~Bourgain and L.~Tzafriri.
\newblock On a problem of {K}adison and {S}inger.
\newblock {\em J. reine angew. Math.}, 420:1--43, 1991.

\bibitem[Buc01]{Buc01:Operator-Khintchine}
A.~Buchholz.
\newblock Operator {K}hintchine inequality in non-commutative probability.
\newblock {\em Math. Annalen}, 319:1--16, 2001.

\bibitem[CR06]{CR06:Quantitative-Robust}
E.~J. Cand{\`e}s and J.~Romberg.
\newblock Quantitative robust uncertainty principles and optimally sparse
  decompositions.
\newblock {\em Foundations of Comput. Math}, 2006.
\newblock To appear.

\bibitem[CT06]{CT06:Kadison-Singer-Problem}
P.~G. Casazza and J.~C. Tremain.
\newblock The {K}adison--{S}inger problem in mathematics and engineering.
\newblock {\em Proc. Natl. Acad. Sci.}, 103(7):2032--2039, Feb. 2006.

\bibitem[Haa82]{Hag82:Best-Constants}
U.~Haagerup.
\newblock The best constants in the {K}hintchine inequality.
\newblock {\em Studia Math.}, 70:231--283, 1982.

\bibitem[JS68]{JS68:Monotone-Convergence}
K.~Jogdeo and S.~M. Samuels.
\newblock Monotone convergence of binomial probabilities and generalization of
  {R}amanujan's equation.
\newblock {\em Ann. Math. Stat.}, 39:1191--1195, 1968.

\bibitem[LP86]{L-P86:Inegalites-Khintchine}
F.~Lust-Picquard.
\newblock In{\'e}galit{\'e}s de {K}hintchine dans {$C_p$ $(1 < p < \infty)$}.
\newblock {\em Comptes Rendus Acad. Sci. Paris, S{\'e}rie I}, 303(7):289--292,
  1986.

\bibitem[LT91]{LT91:Probability-Banach}
M.~Ledoux and M.~Talagrand.
\newblock {\em Probability in Banach Spaces: Isoperimetry and Processes}.
\newblock Springer, 1991.

\bibitem[Rud99]{Rud99:Random-Vectors}
M.~Rudelson.
\newblock Random vectors in the isotropic position.
\newblock {\em J. Functional Anal.}, 164:60--72, 1999.

\bibitem[RV07]{RV07:Sampling-Large}
M.~Rudelson and R.~Vershynin.
\newblock Sampling from large matrices: {A}n approach through geometric
  functional analysis.
\newblock To appear, \textit{J. Assoc. Comput. Mach.}, 2007.

\bibitem[Tim63]{Tim63:Theory-Approximation}
A.~F. Timan.
\newblock {\em Theory of approximation of functions of a real variable}.
\newblock Pergamon, 1963.

\end{thebibliography}

\end{document}